\documentclass[12pt,a4paper]{article}
\usepackage{amsmath,amssymb,latexsym}

\usepackage{graphicx}
\usepackage{amsthm}

\usepackage{wasysym, amsfonts}

\newcommand{\X}{\mathfrak X}
\newcommand{\romb}{\diamondsuit}
\newcommand{\Log}{{Log}}
\newcommand{\MLog}{{\mathbf{ML}}}
\newcommand{\logic}[1]{\mathsf{#1}}

\newcommand{\lK}{\logic{K}}

\newcommand{\Lo}{\mathsf{L}}
\newcommand{\QLo}{\mathsf{QL}}
\newcommand{\rest}[1]{\bigl. #1 \bigr|}

\newcommand{\pmor}{\twoheadrightarrow}

\newcommand{\set}[1]{\left\{#1\right\}}
\newcommand{\Set}[1]{\bigl\{#1\bigr\}}
\newcommand{\setdef}[2][x]{\set{#1\,\left |\,#2 \right .}}
\newcommand{\Setdef}[2][x]{\Set{\bigl.#1\,\bigr|\,#2 }}

\newcommand{\Real}{\mathbb R}
\newcommand{\QQ}{\mathbb Q}
\newcommand{\Natr}{\mathbb N}

\newcommand{\eps}{\varepsilon}

\newcommand{\Y}{\mathcal Y}

\newcommand{\Nf}{\mathcal{N_\omega}}
\newcommand{\mN}{\mathcal{N}}

\newcommand{\len}{\mathop{len}}

\newcommand{\st}{\mathop{st}}

\newcommand{\MF}{\mathcal{MF}}

\newcommand{\pF}{\mathbb{F}}
\newcommand{\pM}{\mathbb{M}}

\newcommand{\pX}{\mathbb{X}}

\title%[N-completeness for quantified modal logic]
{Neighbourhood completeness for quantified pretransitive modal logics}

\author%[A. Kudinov]
{ANDREY KUDINOV}
\newtheorem{theorem}{Theorem}[section]
\newtheorem{thm}[theorem]{Theorem}

\newtheorem{prop}[theorem]{Proposition}

\newtheorem{cor}[theorem]{Corollary}

\newtheorem{lem}[theorem]{Lemma}

\theoremstyle{definition}

\newtheorem{definition}[theorem]{Definition}

\theoremstyle{remark}

\newtheorem*{remark*}{Remark}

\begin{document}
	
	\maketitle
	
	%\begin{frontmatter}

	\begin{abstract}
		We consider quantified pretransitive Horn modal logic. It is known that such logics are complete with respect to predicate Kripke frames with expanding domains. In this paper we prove that they are also complete with respect to neighbourhood frames with constant domains.
	\end{abstract}
	
	%\begin{keyword}
	%  modal logic, predicate logic, neighbourhood semantics, pretransitive logics, Horn axioms.
	%\end{keyword}
	%\end{frontmatter}
	
	\section{Introduction}
	
	In this paper we consider neighbourhood semantics for predicate modal logic.   Neighbourhood semantics for predicate modal logic is drastically
	understudied. There are only a few completeness results for particular logics and no general completeness results. 
	In this paper we aim to take a first step and to prove completeness for a simple jet infinite class of logics. 
		
	Probably, the most straightforward generalization of Kripke semantics to the predicate case is the Kripke semantics with constant domains.
	But in any Kripke frame with a constant domain the Barcan formula $ \forall x \left (\Box P(x)\right ) \to \Box \forall x P(x) $ is valid. 
	At the same time, the Barcan formula is not derivable in predicate modal logics in general. Instead, Kripke frames with expanding domains were considered, and several completeness results were proved (see \cite{gabbay2009quantification}). 
	But still, many simple predicate modal logics are Kripke incomplete. Different authors suggested several generalizations, including neighbourhood semantics, (for an overview see \cite[Part II]{gabbay2009quantification}). 
	
	In the case of transitive reflexive logics, neighbourhood semantics is equivalent to the topological semantics. \cite{rasiowa1963metamathematics} showed that $ \logic{QS4} $ is complete with respect to topological spaces with constant domain. 
	In a recent paper, \cite{kremer2014quantified} has proved that $ \logic{QS4} $ is complete with respect to the set of rational numbers $ \QQ $.  
	\cite{arlo2006first} have proved that the minimal normal modal predicate logic $ \logic{QK} $ is complete with respect to neighbourhood frames with constant domains. The methods for proving completeness were different in the last two papers. In \cite{arlo2006first} authors used a neighbourhood canonical model construction, whereas  \cite{kremer2014quantified} used the Kripke semantics as an intermediate step and constructed p-morphisms from $ \QQ $ onto a Kripke frame with expanding domains.
	
	The intuition behind the Kremer's construction is that a topological (or, more generally, neighbourhood) frame with a constant domain is very similar to a product of topological space and an $\logic{S5}$-frame, where the $\logic{S5}$-frame plays the role of the domain. 
	%Here we understand the notion ``product'' as in \cite{benthem:MultimodalLogicsProductsTopologies}, note that it differs from the classical notion of the product in topology, the details will be given further in this paper. 
	In this paper we use the same idea and prove that $ \QLo $ is complete with respect to neighourhood frames with constant domain, where $ \Lo $ is a one-way PTC-logic (see Definition \ref{def:PTClogic}).
	We adopt the methods developed by the author for neighbourhood-Kripke products in \cite{Kudinov_KNproduct17}.

	\section{Propositional modal logic}
	\subsection{Syntax}
	Let $ \mathrm{PROP} $ be a countable set of propositional letters. A \emph{(propositional modal) formula} is defined recursively using the Backus-Naur form as follows:
	$$
	A ::= p\; |\;\bot \; | \; (A \to A) \; | \; \Box_i A,
	$$
	where $p\in \mathrm{PROP}$, and $\Box_i$ is a modal operator ($ i= 1, \ldots, N$).
	Other connectives are introduced as abbreviations: classical connectives are expressed through $\bot$ and $\to$,
	and $\Diamond_i$ is a shorthand for $\lnot\Box_i\lnot$. The set of all modal formulas is denoted by $ \mathcal{ML}_N $.
	%, and in order to specify the modalities used in the language we write them in subscript, for example: $ \mathcal{ML}_{\Box_1} $ or $ \mathcal{ML}_{\Box_2} $. 
	
	\begin{definition} A \emph{normal modal logic\/} (or \emph{a logic}, for short) is
		a set of modal formulas closed under Substitution $\left
		(\frac{A(p)}{A(B)}\right )$, Modus Ponens $\left (\frac{A,\,
			A\to B}{B}\right )$ and Generalization rules $\left
		(\frac{A}{\Box_i A}\right)$, containing all the
		classical tautologies and the normality axioms:
		$$
		\begin {array}{l}
		\Box_i (p\to q)\to (\Box_i p\to \Box_i q).
	\end{array}
	$$
	
	$\logic{K_N}$ denotes \emph{the minimal normal modal logic with $n$ modalities} and $\lK = \logic{K_1}$.
\end{definition}

Let $\logic{L}$ be a logic and $\Gamma$ be a set of formulas, then $\logic{L} +
\Gamma$ denotes the minimal logic containing $\logic{L}$ and $\Gamma$. If
$\Gamma = \set{A}$, then we write $\logic{L} + A$ rather than $\logic{L}+  \{ A \}$.

\subsection{Kripke semantics}
\begin{definition}\label{def:Kripke_model}
	A \emph{Kripke frame} is a tuple $F = (W, R)$, where $ W $ is a non-empty set and $ R \subseteq W \times W $.
	
	For a Kripke frame $F = (W, R)$ we define the \emph{subframe} generated by $w \in W$ as the frame $F^w = (W', R|_{W'})$, where $W' = R^*(w)$ and ${R}|_{W'} = R \cap W' \times W'$. The star $ {}^* $ here stands for the reflexive and transitive closure of a relation. A frame $F$ is \emph{rooted} if $F=F^w$ for some $w$.
	
	A frame $ F $ with a \emph{valuation} $ V:PROP \to 2^W $ is \emph{a model} $M = (F, V)$.  We say that model $ M $ is \emph{based} on frame $ F $.
	
	The truth of a formula in a model $ M $ at a point $ x\in W $ is defined, as usual, by induction on the length of the formula:
	\begin{align*}
		M, x &\not\models \bot; &\\
		M, x &\models p &\iff& x \in V(p);\\
		M, x &\models A \to B &\iff& M, x \not\models A \hbox{ or } M, x \models B;\\
		M, x &\models \Box A &\iff& \forall y\; (xR y \Rightarrow M, y \models A).
	\end{align*} 
	
	%Sometimes it is useful to define the truth notion in a model in terms of extending the valuation to the set of all formulas in the following way:
	
	A formula $ A $ is \emph{true in a (Kripke) model} $M$ if $\forall x\in W (M, x \models A)$ (notation $M \models A$).
	
	A formula $ A $ is \emph{valid on a (Kripke) frame} $F$ if for any valuation $ V $ $ (F, V) \models A$ (notation $F \models A$). 
	
	A formula is \emph{valid on a (Kripke) frame $ F $ at a point $ x $}  if it is true in all models based on $F$ at point $ x $ (notation $F, x \models A$). That is if $ \forall V (F,V, x \models A) $.

	We write $F \models \Lo$ if $F\models A$ for all $A \in \Lo$.
	%The \emph{logic} of a class of Kripke frames $\mathcal{C}$ is $\Log(\mathcal{C}) = \setdef[A]{F \models A \hbox{ for all } F \in \mathcal{C}}$.
	%For a logic $\Lo$ we also define $\mathcal{V}(\Lo) =\setdef[F]{\hbox{$F$ is a Kripke frame and } F\models \Lo}$. Note that if there are no $F$ such that $F \models \Lo$, then $\mathcal{V}(\Lo) = \varnothing$.
\end{definition}

We define the logic of a class of Kripke frames $\mathcal{F}$ as 
$$\Log(\mathcal{F}) = \setdef[A]{F \models A \hbox{ for all } F \in \mathcal{F}}.$$

For a single frame we write $\Log(F)$ instead of $\Log(\set{F})$.

%For a logic $\logic{L}$ we also define 
%$$
%\kV(\logic{L}) =\setdef[F]{\hbox{$F$ is an Kripke frame and } F\models \logic{L}}.
%$$
%Note that if there are no $F$ such that $F \models \logic{L}$, then $\kV(\logic{L}) = \varnothing$.

%\begin{remark*}
We sometimes write $w\in F$ as a shorthand for ``$w \in W$ and $F = (W, R)$''.
%\end{remark*}

\begin{definition} \label{def:p_morphism}
	Let $F = (W, R)$ and $G = (U, S)$ be Kripke frames. A function $f: W \to U$ is a \emph{p-morphism} (Notation: $f: F \pmor G$) if 
	\begin{enumerate}
		\item $f$ is surjective;
		\item {}\textbf{[monotonisity]} for any $w, v\in W$ $ wR v $ implies $ f(w) S f(v) $;
		\item {}\textbf{[lifting]} for any $w\in W$ and $v' \in U$ such that $ f(w)S v' $ there exists $v \in W$ such that $ wRv $ and $ f(v)=v' $.
	\end{enumerate}
\end{definition}

The following p-morphism lemma is well-known (see \cite[Proposition 2.14]{blackburn_modal_2002})

\begin{lem}\label{lem:pmorhism4K-frame}
	Let $f: F \pmor G$ and $V'$ be a valuation on $G$. We define a valuation $ V $ on $ F $, such that $V(p) = f^{-1}(V'(p))$. Then, for any $ x\in F $ and formula $ A $
	\[
	F, V, x \models A \iff G, V', f(x) \models A.
	\] 
\end{lem}

%The proof is by standard induction on the length of the formula $A$. 
%The following is a straightforward corollary. 

%\begin{corollary}\label{cor:pmorphKframes}
%	If $f: F \pmor F$, then $\Log(F) \subseteq \Log(G)$.
%\end{corollary}	

\begin{definition}
	A formula $A$ is \emph{closed} if it contains no variables.
\end{definition}

\begin{lem}
	For any closed modal formula $ A $ and a p-morphism of Kripke frames $f:F \pmor G$ 
	\[ 
	F,x \models A \iff G, f(x) \models A.
	\]
\end{lem}

This follows from Lemma \ref{lem:pmorhism4K-frame} since the truth of a closed formula does not depend on the valuation.

%\begin{definition}    
%	Following \cite{GSh_Product_1998}, we define \emph{a universal strict Horn sentence} as a first-order closed formula of the form
%	\[ 
%	\forall x \forall y\forall z_1 \ldots \forall z_n \bigl(\phi(x,y, z_1, \ldots, z_n) \to \psi(x,y)\bigr),
%	\]
%	where $ \phi(x,y, z_1, \ldots, z_n) $ is quantifier-free positive (i.e.{}, it is built from atomic formulas by using $ \land $ and $ \lor $) and $ \psi(x,y) $ is an atomic formula in the signature $ \Omega = \left\langle R^{(2)}_1, \ldots, R^{(2)}_m\right \rangle $, where $ R^{(2)}_i $ is the propositional letter that corresponds to the relation $ R_i $.
%\end{definition}

%	 By a \emph{closed sentence} we understand the standard translation of a closed modal formula\footnote{This definition is rather controversial since in 1st order logic sentence is already a closed formula, but we did not find a better name for this notion.}.

%	 A logic $ \logic{L} $ is \emph{Horn axiomatizable} if the class of all its frames is the class of all 1st-order models of a set of Horn and closed sentences.

We define formula $ \Box^k A$ by induction:
\begin{itemize}
	\item $ \Box^0 A = A $, here $ \eps $ is the empty sequence;
	\item $ \Box^{k+1} A = \Box^k \Box A $.
\end{itemize} 

\begin{definition}
	Logic $ \logic{L} $ is called \emph{pretransitive} if for some $ k>0 $ formula 
	$ p \land \Box p \land \ldots \land \Box^{k} p \to \Box^{k+1} p  $ is in $ \logic{L} $.
\end{definition}

Sometimes such logics are called \emph{weakly transitive}  \cite[Section 3.4]{kracht1999tools}.

\begin{definition}\label{def:PTClogic}
	A logic $ \Lo $ is a \emph{one-way PTC-logic}\footnote{In \cite{Kudinov_KNproduct17} such logics are called HTC-logics}  if it can be axiomatized by closed formulas and formulas of the type $ \Box p \to \Box^k p $, where $ k\ge 0 $. These formulas correspond to universal strict Horn sentences (see \cite{GSh_Product_1998}). 
\end{definition}

\begin{remark*}
	It is easy to see that a one-way PTC-logic is pretransitive.
\end{remark*}

The following well-known and well-studied logics are one-way PTC-logic:
\begin{align*}
	\logic{D} &= \logic{K} + \lnot \Box \bot, & \logic{T} &= \logic{K} + \Box p \to p,\\ 
		\logic{K4} &= \logic{K} +  \Box p \to \Box\Box p, & \logic{S4} &= \logic{T} + \Box p \to \Box\Box p,\\
		\logic{D4} &= \logic{K4} + \lnot \Box \bot.
\end{align*}
Whereas logic $ \logic{S5} = \logic{S4} + p \to \Box \romb p$ is not a one-way PTC-logic.

\begin{definition}    
	Following \cite{GSh_Product_1998}, we define \emph{a universal strict Horn sentence} as a first-order closed formula
	\[ 
	\forall x \forall y\forall z_1 \ldots \forall z_n \bigl(\phi(x,y, z_1, \ldots, z_n) \to \psi(x,y)\bigr),
	\]
	where $ \phi(x,y, z_1, \ldots, z_n) $ is quantifier-free
	positive (i.e.{}, it is built from atomic formulas using $ \land $ and $ \lor $) and $ \psi(x,y) $ is an atomic formula in the signature $ \Omega = \left\langle R^{(2)}_1, \ldots, R^{(2)}_m\right \rangle $, where $ R^{(2)}_i $ is a 2-ary propositional letter. In our case $ m=1 $.
\end{definition}

Let $ \Gamma $ be a set of universal strict Horn sentences and $F$ be a Kripke frame.  We define the $ \Gamma $-closure of $ F $ by $ F^{\Gamma} $. It is the minimal (in terms of inclusion of relations) frame such that all sentences from $ \Gamma $ are valid in it. Such a frame exists due to

\begin{lem}[{\cite[Prop 7.9]{GSh_Product_1998}}]\label{lem:HornClosure}
	For any Kripke frame $ F=(W, R_1, \ldots, R_N) $ and a set of universal strict Horn sentences $\Gamma$, there exists $ F^{\Gamma} = (W, R^{\Gamma}) $ such that
	\begin{itemize}
		\item $ R \subseteq R^{\Gamma} $;
		\item $ F^{\Gamma} \models \Gamma $;
		\item if $ G \models \Gamma $ and $ f:F \pmor G $ then $ f: F^{\Gamma} \pmor G$.
	\end{itemize}
\end{lem}

The minimality of $ F^\Gamma $ follows from the proof.

\begin{definition}
	Let $ F = (W, R)$ be a Kripke frame. A path in $ F $ is a tuple $ w_0 R w_1 \ldots R w_m $, where for all $ i\in\set{0,\ldots, m-1} w_i R w_{i+1}$.
\end{definition}

\begin{definition}\label{def:unraveling}
Let $ F = (W, R)$ be a rooted frame with root $ w_0 $.
A path in $ F $ is \emph{rooted} if it starts with $ w_0 $. 
Let $ W^{\sharp} $ be the set of all rooted paths in $ F $.
	For any rooted path $ \alpha = w_0 R w_1 \ldots R w_m $ we define 
\begin{align*}
	\pi (\alpha) &= w_m;\\
	\alpha R^{\sharp} \beta &\iff \beta = \alpha w_{m+1} \hbox{ and } \pi(\alpha) R \pi(\beta).
\end{align*}	
Frame $ 	F^{\sharp} = (W^{\sharp}, R^{\sharp}) $ is the \emph{unravelling} of $F$.
\end{definition}

\begin{lem}\label{lem:unrev_pmorph}
	Map $ \pi $ is a p-morphism: $ \pi: F^{\sharp} \pmor F $.
\end{lem}

The proof is straightforward (c.f. \cite[Lemma 4.52]{blackburn_modal_2002}).

%Further in this paper we will apply Lemma \ref{lem:HornClosure} to the p-morphism from Lemma \ref{lem:unrev_pmorph}.

\subsection{Neighborhood semantics}
We will only give the definitions we need for this paper, more information on this topic and general definitions can be found in \cite{Segerberg1971,Chellas1980} and a more recent book \cite{pacuit2017neighborhood}.

\begin{definition}\label{def:filter}
	Let $X$ be a non-empty set of \emph{points}, then $\mathcal{F} \subseteq 2^X$ is a \emph{filter} on $X$ if
	\begin{enumerate}
		\item $X \in \mathcal{F}$;
		\item if $U_1,\; U_2 \in \mathcal{F}$, then $U_1 \cap U_2 \in \mathcal{F}$;
		\item if $U_1 \in F$ and $U_1 \subseteq U_2$, then $U_2 \in \mathcal{F}$.
	\end{enumerate}
	It is usually required that $\varnothing \notin \mathcal{F}$ ($\mathcal{F}$ is a proper filter), but we will not require it in this paper.
\end{definition}

\begin{definition}\label{def:filterbase}
	For $X \ne \varnothing$ a set of subsets $ \mathcal{B} \subseteq 2^X$ is a \emph{filter base} if 
	\begin{enumerate}
		\item $ \mathcal{B} \ne \varnothing $;
		\item for any $U_1,\; U_2 \in \mathcal{B}$  $ \exists U_3 \in \mathcal{B} \left ( U_3 \subseteq U_1 \cap U_2 \right )$.
	\end{enumerate}
	Given a filter base $ \mathcal{B} $, the filter \emph{generated} by $ \mathcal{B} $ is defined as the minimal filter containing $ \mathcal{B} $. It is the family of all supersets of sets from $ \mathcal{B} $. 
\end{definition}

\begin{definition}
	A \emph{(normal) neighbourhood frame} (an \emph{n-frame} for short) is a tuple $\X = (X, \tau)$, where $X$ is a nonempty set and $\tau: X \to 2^{2^X}$ such that $\tau(x)$ is a filter\footnote{Usually neighbourhood semantics is used for non-normal logics, and in the most general case there are no restrictions on the neighbourhood function, but here we will consider only normal modal logics.} on $X$ for any $x$. The function $\tau$ is called the \emph{neighbourhood function} of $\X$, and elements of $\tau(x)$ are called \emph{neighbourhoods of $x$}.
	
	\emph{A neighbourhood model (or n-model)} is a pair $(\X, V)$, where $\X = (X, \tau)$ is an n-frame and $V: PROP \to 2^X$ is a \emph{valuation}.  We say that model $ (\X, V) $ is based on $ \X $.
\end{definition}

%\begin{remark}
%	Many authors consider neighbourhood semantics for non-normal and even non-monotone logics. In this case a set of neighbourhoods can be an arbitrary set of sets. Other authors consider  \emph{monotone} neighbourhood frames and require only item 3 from Definition \ref{def:filter} (a set of neighbourhoods is closed under supersets). 
%\end{remark}

\begin{definition}
	The \emph{truth of a formula in a neighbourhood model} is defined by induction. For the variables and Boolean connectives the definition is the same as for Kripke model (Def. \ref{def:Kripke_model}). For modalities the definition is the following:
	\[
	M, x \models \Box A \iff \exists U \in \tau(x)\; \forall y \in U (M,y \models A).
	\]
	A formula is true in an n-model $M$ if it is valid at all points of $M$ (notation $M \models A$).
	A formula is valid on an n-frame $\X$ if it is true in all models based on $\X$ (notation $\X \models A$). 
	We write $\X \models \logic{L}$ if for any $A \in \logic{L}$, $\X\models A$.
	We define the logic of a class of n-frames $\mathcal{C}$ as $\Log(\mathcal{C}) = \setdef[A]{\X \models A \hbox{ for all } \X \in \mathcal{C}}$ and $ \Log(\X) = \Log(\set{\X}) $.
	%For a logic $\logic{L}$ we also define $\nV(\logic{L}) =\setdef[\X]{\hbox{$\X$ is an n-frame and } \X\models \logic{L}}$. Note that if there are no $\X$ such that $\X \models \logic{L}$, then $\nV(\logic{L}) = \varnothing$.
\end{definition}

Given a Kripke frame one can construct an equivalent n-frame:
\begin{definition}
	Let $F = (W, R)$ be a Kripke frame. Then $\mN(F) = (W, \tau)$ is a n-frame, such that
	\[
	\tau (w) = \setdef[U]{ R(w) \subseteq U \subseteq W},
	\] 
	where $ R(w) = \setdef[u]{wRu} $.
\end{definition}

Frames $ F $ and $ \mN(F) $ are equivalent in the following sense 
\begin{lem}\label{lem:n-frame_from_kframe}
	Let $F = (W, R)$ be a Kripke frame. Then
	\[
	\Log(\mN(F)) = \Log(F).
	\]
\end{lem}

The proof is straightforward (see \cite{Chellas1980}). 

This lemma shows that neighbourhood semantics is a generalization of Kripke semantics. 

\begin{definition} \label{def:bounded_morphism}
	Let $\X = (X, \tau)$ and $\Y = (Y, \sigma)$ be neighbourhood frames. Then function $f: X \to Y$ is a \emph{p-morphism} (notation $f: \X \pmor \Y$) if 
	\begin{enumerate}
		\item $f$ is surjective;
		\item \textbf{[zig]} for any $x\in X$ and $U \in \tau(x)$, we have $f(U) \in \sigma(f(x))$;
		\item \textbf{[zag]} for any $x\in X$ and $V \in \sigma(f(x))$, we have $f^{-1}(V) \in \tau(x)$.
	\end{enumerate}
\end{definition}

%\begin{remark}
%	A neighbourhood frame is a generalization of a topological space. In paticular for a topological space $ (X, T) $ we can define neighbourhood function $ \tau_T (x) = \setdef[U]{\exists U' \in T (x\in U' \land U' \subseteq U)} $. Then $ (X, \tau_T) $ is a n-frame 
%\end{remark}

%\begin{remark}
%	According to Lemma \ref{lem:n-frame_from_kframe}, a Kripke frame is a special case of a neighbourhood frame. It is easy to check that for any two Kripke frames $F$ and $G$ a function $f$ is a p-morphism (Definition \ref{def:p_morphism}) from $F$ to $G$ iff $f$ is a p-morphism (Definition \ref{def:bounded_morphism}) from $\mN(F)$ to $\mN(G)$. So, a p-morphism for n-frames is a natural generalization of the notion of p-morphism for Kripke frames. This is why we use the same name for these two formally different notions.
%\end{remark}

\begin{lem}\label{lem:pmorhism4n-frame}
	Let $\X = (X, \tau)$, $\Y = (Y, \sigma)$ be n-frames and $f: \X \pmor \Y$. Let $V'$ be a valuation on $\Y$. We define valuation $V(p) = f^{-1}(V'(p))$. Then
	\[
	\X, V, x \models A \iff \Y, V', f(x) \models A.
	\] 
\end{lem}

The proof is by induction on the length of formula $A$. The following is a straightforward corollary. 

\begin{cor}\label{cor:pmorph}
	If $f: \X \pmor \Y$, then $\Log(\X) \subseteq \Log(\Y)$.
\end{cor}

\section{Predicate modal logic}
%All results in this paper can be formulated from multimodal case and the proofs should similar but the details should be checked carefully.
%Here however we consider unimodal case. 
Following \cite{gabbay2009quantification} we define predicate modal formulas and logics as follows.

\begin{definition}
	Let $ \mathop{Var} $ be a countably infinite set of \emph{(individual) variables} and $ PL^n = \setdef[P^n_i]{i\ge 0} $ be a fixed set of \emph{n-ary predicate letters} ($ n\ge 0 $). 0-ary predicate letters we will call propositional letters.
	
	Modal predicate \emph{formulas} are defined inductively as follows:
	
	\begin{itemize}
		\item $ \bot $ is a formula;
		\item $ P^0_i $ is a formula;
		\item if $ x_1, \ldots, x_k \in Var$ then $ P^k_i(x_1, \ldots, x_k) $ is a formula;
		\item if $ A $ and $ B $ are formulas then $ (A \to B) $ is a formula;
		\item if $ A $ is a formula then $ \Box A $ is a formula;% ($ i \le N $);
		\item if $ A $ is a formula and $ x \in Var $ then $ \forall x A $ is a formula.
	\end{itemize}
	
	All other connectives $ \land, \lor, \lnot, \exists, \romb $ are expressed as usual.
	
	The set of all modal predicate formulas is denoted by $ \MF $.
\end{definition}

%Let $ \Omega = \langle Pred, Func \rangle$ be a signature, where $Pred = \bigcup_{n\ge 1} Pred_n \ne \varnothing $, $ Pred_n $ is the set of n-ary predicate letters, and $ Func = \bigcup_{n \ge 0} Func_n$, $ Func_n $ is the set of n-ary function letters. Elements of $ Func_0 $ we will call constants. A formula of predicate modal logic is constructed using infinite set of variables $ Var $, symbols from $ \Omega $, logical connectives: $ \to $, $ \bot $ (constant ``False'') and quantifier $ \forall $, modal operators $ \set{\Box_i}_{i \le N} $ and brackets in the following way.

%\emph{Terms} are defined inductively as follows:
%\begin{itemize}
%	\item a variable (element of $ Var $) is a term;
%	\item if $ t_1, \ldots, t_k $ are terms and $ f \in Func_k $ then $ f(t_1, \ldots, t_k) $ is a term (note that it is possible that $ k=0 $, in this case $ f $ by itself is a term).
%\end{itemize}

\begin{definition}
	If $ \Lo $ is a modal logic then $ \QLo $ is the minimal set of predicate modal formulas such that
	\begin{itemize}
		\item $ \QLo $ includes all formulas from $ \Lo $ where propositional variables are replaced by corresponding propositional letters.
		\item  for a 1-ary predicate $ P $ and a propositional letter $ Q $,  $ \QLo $ includes formulas  
		\begin{itemize}
			\item $ \forall x P(x) \to P(y) $, 
			\item $ \forall x (Q \to P(x)) \to (Q \to \forall x P(x)) $;
			\item $\forall x ( P(x) \to Q ) \to (\exists x P(x) \to Q) $.
		\end{itemize}
		\item $\QLo$ is closed under Modus Ponens $ \left ( \frac{A,\ A \to B}{B}\right ) $, necessitation $ \left ( \frac{A}{\Box A}\right ) $, universal generalization rules $ \left ( \frac{A}{\forall x A}\right ) $, and under $ \MF $-substitutions \footnote{The idea is that we avoid variable collisions. For the sake of simplicity one can think that we rename all bound variable before substitution. Also we will assume that two formulas are congruent if they are the same up to renaming the bound variables.} (a proper definition of these substitutions can be found in \cite{gabbay2009quantification}). 
	\end{itemize}
\end{definition}

\begin{lem}[see \cite{gabbay2009quantification}]
	For any modal logic $ \Lo $ the quantified modal logic $ \QLo $ includes formula $ \Box \forall x P(x) \to \forall x \Box P(x)$ (the converse Barcan formula).
\end{lem}

Note that $ \QLo $ in general does not include the Barcan formula: 
$$\forall x \Box P(x) \to \Box \forall x  P(x).$$

\section{Semantics for predicate modal logic}

\begin{definition}
	A system of expanding domains over a Kripke frame $ F = (W, R) $ is a family of sets $ D = (D_u)_{u \in W} $ such that 
	\[  
	\forall u, v \in W (u R v \Rightarrow D_u \subseteq D_v).
	\]
	
	A \emph{predicate Kripke frame} (with expanding domains) is a pair $ \pF = (F, D) $, where $ F $ is a Kripke frame and $ D $ is a \emph{system of expanding domains} over $ F $. 
\end{definition}

\begin{definition}
	%For a Kripke frame $F = (W, R_1, \ldots, R_N)$ we define the \emph{subframe} generated by $w \in W$ as the frame $F^w = (W', R_1|_{W'}, \ldots, R_N|_{W'})$, where $W' = (R_1 \cup \ldots, \cup R_N)^*(w)$ and ${R_i}|_{W'} = R_i \cap W' \times W'$. A frame $F$ is \emph{rooted} if $F=F^w$ for some $w$.	
	A \emph{subframe} of a predicate frame $ \pF = (F, D) $ generated by $w \in F$ is a predicate frame $\pF^w = (F^w, D')  $ such that $ F^w = (W', R') $ is the subframe of $ F $ generated by $ w $ (see Def. \ref{def:Kripke_model}) and $ D' = D|_{W'} = \setdef[D_u]{u \in W'}$.
\end{definition}

\begin{definition}
	A \emph{valuation} $ \xi $ on a predicate frame $ \pF $ is a function sending every predicate letter $ P^m_k $ to a family of m-ary predicates on the domains:
	\[ 
	\xi(P^m_k) = (\xi_u(P^m_k))_{u\in W}, \hbox{ where } \xi_u(P^m_k) \subseteq D_{u}^{m}.
	\]
	A frame $ \pF $ with a valuation $\xi$ is \emph{a (predicate) model} $\pM = (\pF, \xi)$.
	Note that for $ m=0 $  $ \xi_u(P^0_k) $ is an 0-ary predicate so it is ether true or false. It is very similar to propositional variables.
	
	The truth of a closed formula (formula without free variables) in a model $ \pM $ at a point $ u\in W $ is defined, by induction on the length of the formula, but first we enrich our language with constants from set $ \bigcup_{u\in W} D_u $:
	\begin{align*}
		\pM, u &\not\models \bot; &\\
		\pM, u &\models P^0_i &\iff& \xi_u(P^0_i) \hbox{ is true};\\
		\pM, u &\models P^m_i(a_1, \ldots, a_m) &\iff& (a_1, \ldots, a_m) \in \xi_u(P^m_i);\\
		\pM, u &\models A \to B &\iff& \pM, u \not\models A \hbox{ or } \pM, u \models B;\\
		\pM, u &\models \Box A &\iff& \forall v\; \left (uR v \Rightarrow \pM, v \models A\right );\\
		\pM, u &\models \forall x A(x)&\iff &\forall a \in D_u \left (\pM, u \models A(a)\right ).
	\end{align*} 
\end{definition}

Let $ A(x_1, \ldots, x_k) $ be a formula with free variables such that it does not have free variables other then $ x_1, \ldots, x_k $. Then the universal closure of $ A(x_1, \ldots, x_k) $ is the following closed formula
\[ 
\bar\forall A = \forall x_1 \ldots \forall x_k\; A(x_1, \ldots, x_k)
\]

A formula is \emph{true in a (Kripke) model} $\pM$ if its universal closure is true at all points of $\pM$ (notation $\pM \models A$).
A formula is \emph{valid on a (Kripke) frame} $\pF$ if it is true in all models based on $\pF$ (notation $\pF \models A$). 
We write $\pF \models \Gamma$ if, for any $A \in \Gamma$, $\pF\models A$.
The \emph{logic} of a class of Kripke frames $\mathcal{C}$ is $\MLog(\mathcal{C}) = \setdef[A]{\pF \models A \hbox{ for all } \pF \in \mathcal{C}}$.
Logic $ \Lo $ is \emph{complete} with respect to Kripke semantics with expanding domains if there is a class of Kripke frames with expanding domains $ \mathcal{C} $ such that $ \MLog(\mathcal{C}) = \Lo $.
%For a logic $\Lo$ we also define $\mathcal{V}(\Lo) =\setdef[F]{\hbox{$F$ is a Kripke frame and } F\models \Lo}$. Note that if there are no $F$ such that $F \models \Lo$, then $\mathcal{V}(\Lo) = \varnothing$.

%\begin{lem}
%	$f: \cF \pmor \cF^0$.
%\end{lem}
%\begin{proof}
%\textbf{Monotonicity.}  Assume that $ x R_i y $ then $ \setdef[A]{\Box_i A \in x} \subset y$ and
%\[ 
%\setdef[A]{\Box_i A \in \bar x} = \setdef[A]{\Box_i A \in x} \cap \Upsilon \subset y \cap \Upsilon = \bar y.
%\]
%So it follows that $ \bar x R'_i \bar y$.

%\textbf{Lifting.} Assume that $ \bar x R'_i \bar y $ consider set
%\[ 
%\Gamma = \bar y \cup \setdef[A]{\Box A \in x}.
%\]

%This set is consistent, indeed 
%\end{proof}

%Here we describe a construction of continuum unravelling. It is similar to the construction in \cite[Lemma 4.9]{GSh_Product_1998}.

%\begin{definition}
%Let $ F=(W, R) = F^{w_0}$ be a rooted Kripke frame, $ S $ be a non-empty set and $ x_0 \in S $ be a point in it (the starting point). Then
%\begin{align*}
%F\cdot S &= (W\times S, R \cdot S);\\
%(w,x) R \cdot S (v,y) &\iff wRv;\\
%F_S = (F\cdot S)^{(w_0, x_0)} &= (W_S, R_S) \hbox{ --- a rooted subframe}.
%\end{align*}
%$ F_S $ is called the \emph{thickening} of $ F $ by $ S $.
%\end{definition}

\begin{thm}[{\cite[Theorem 6.1.29]{gabbay2009quantification}}]\label{thm:PTC_ncompleteness}
	The quantification $ \logic{QL} $ of any one-way PTC-logic $ \logic{L} $ is complete with respect to Kripke semantics with expanding domains.
\end{thm}

%This theorem is proved using the canonical model. Our language is countable, so the canonical model is continuum and the cardinality of domains in it is no greater then continuum. 

\begin{definition}
	Let $ \pF = (W, R, D)$ and $ \pF' = (W', R', D')$ be a Kripke frames with expanding domains, $ D= \bigcup_{w\in W} D_w $, $ D' = \bigcup_{w\in W} D'_w $.
	A \emph{p-morphism} from  $ \pF$ to $ \pF' $ is a pair $ (\varphi_0, \varphi_1) $, such that:
	\begin{enumerate}
		\item $ \varphi_0 : (W, R) \pmor (W', R') $;
		\item $ \varphi_1 = (\varphi_{1w})_{w\in W} $ is a family of surjective functions: $$ \varphi_{1w}:D_{w} \to D_{\varphi_0(w)};$$
		\item if $ w R w'$ then $\forall d \in D_w  \bigl(\varphi_{1w} (d) = \varphi_{1w'}(d)\bigr) $.
	\end{enumerate}
	Notation: $ (\varphi_0, \varphi_1): \pF \pmor \pF' $.
	We write $ \pF \pmor \pF' $ if there exists a p-morphism from $ \pF $ to $ \pF' $.	
\end{definition}

\begin{lem}[\cite{gabbay2009quantification}, Prop. 3.3.11 \& 3.3.12]\label{lem:pmorphismlemma4kF2kF}
	Let $ (\varphi_0, \varphi_1): \pF \pmor \pF' $ and $ \xi' $ be a valuation on $ \pF' $. We define valuation $ \xi$ on $ \pF $ in the following way 
	$$ 
	(a_1, \ldots, a_m) \in \xi_u(P^m_k)  \iff \left (\varphi_{1u}(a_1), \ldots, \varphi_{1u}(a_m)\right ) \in \xi'_{\varphi_0(u)}(P^m_k).
	$$
	Then for any $ u\in W $ and any formula $ A $
	\[ 
	\pF, \xi, u \models \bar \forall A \iff \pF', \xi', \varphi_0(u) \models \bar \forall A.
	\]
\end{lem}

\begin{definition}
	A predicate neighbourhood frame with constant domain is a pair $ \pX = (\X, D^*) $, such that $ \X = (X, \tau)$ is a neighbourhood frame and $ D^* $ is a nonempty set.
	
	A \emph{valuation} $\theta $ on $ \pX $ is a function sending every predicate letter $ P^m_k $ to a family of m-ary predicates on $ D^* $:
	\[ 
	\theta(P^m_k) = (\theta_u(P^m_k))_{u\in W}, \hbox{ where } \theta_u(P^m_k) \subseteq (D^*)^m.
	\]
	A \emph{neighbourhood model} on $ \pX $ is a pair $ \pM = (\pX, \theta) $.
	
	The truth of a closed formula in a model $ \pM $ at a point $ x\in X $ is defined similar to Kripke models, by induction on the length of the formula. We also enrich our language with constants from set $D^*$. 
	\begin{align*}
		\pM, x &\not\models \bot; &\\
		\pM, x &\models P^0_i &\iff& \xi_x(P^0_i) \hbox{ is true};\\
		\pM, x &\models P^m_i(a_1, \ldots, a_m) &\iff& (a_1, \ldots, a_m) \in \xi_x(P^m_i);\\
		\pM, x &\models A \to B &\iff& \pM, x \not\models A \hbox{ or } \pM, x \models B;\\
		\pM, x &\models \Box A &\iff& \exists U \in \tau(x) \forall y \in U  \left ( \pM, y \models A\right );\\
		\pM, x &\models \forall x A(x)&\iff &\forall a \in D^* \left (\pM, x \models A(a)\right ).
	\end{align*} 
\end{definition}

%Since the composition of p-morphism is a p-morphism
%\begin{lem}
%	The map $ f = \pi \circ p_1 : F^{\sharp}_{\Real} \to F$ is a p-morphism.
%\end{lem}

\begin{definition} \label{def:pmorph4prednframe}
	Let $ \pX = (X, \tau, D^*)$ be a neighbourhood frame with constant domain and $ \pF = (W, R, D)$ be a Kripke frame with expanding domains, $ D= \bigcup_{w\in W} D_w $.
	A \emph{p-morphism} from  $ \pX$ to $ \pF $ is a pair $ (\varphi_0, \varphi_1) $, such that:
	\begin{enumerate}
		\item \label{it:1:def:pmorph4prednframe} $ \varphi_0 : (X, \tau) \pmor \mN(W, R) $;
		\item \label{it:2:def:pmorph4prednframe} $ \varphi_1 = (\varphi_{1x})_{x\in X} $ is a family of surjective functions indexed by points from~$X$: $$ \varphi_{1x}:D^* \to D_{\varphi_0(x)};$$
		\item \label{it:3:def:pmorph4prednframe} $ \forall d \in D^*\,\forall x \in X \,\exists U \in \tau(x)\, \forall y \in U \bigl(\varphi_{1y} (d) = \varphi_{1x}(d)\bigr) $.
	\end{enumerate}
	Notation: $ (\varphi_0, \varphi_1): \pX \pmor \pF $.
	We write $ \pX \pmor \pF $ if there exists a p-morphism from $ \pX $ to $ \pF $.	
\end{definition}

\begin{lem}\label{lem:pmorphismlemma4nF2kF}
	Let $ \pX = (X, \tau, D^*)$ be a neighbourhood frame with constant domain, $ \pF = (W, R, D)$ be a Kripke frame with expanding domains, $ (\varphi_0, \varphi_1): \pX \pmor \pF $ and $ \xi $ be a valuation on $ \pF $. We define valuation $ \theta = (\theta_x)_{x \in X} $ on $ \pX $ in the following way 
	$$ 
	(a_1, \ldots, a_m) \in \theta_x(P^m_k)  \iff \left (\varphi_{1x}(a_1), \ldots, \varphi_{1x}(a_m)\right ) \in \xi_{\varphi_0(x)}(P^m_k).
	$$
	Then for any $ x\in X $ and any formula $ A $
	\[ 
	\pX, \theta, x \models \bar \forall A \iff \pF, \xi, \varphi_0(x) \models \bar \forall A.
	\]
\end{lem}
\begin{proof}
	Let us assume that $A$ has no free variables but can contain constants from set $ D^* $ (for models on $ \pX $) and from set $\bigcup D $ (for models on $ \pF $).  Assume that $ A $ includes only constants $ a_1, \ldots, a_m $, to highlight it, we will write $ A(a_1, \ldots a_m) $.

	We prove the following statement by induction on the length of $ A $:
	\[  
	\forall x\in X \bigl (\pX, \theta, x \models  A(a_1, \ldots a_m) \iff \pF, \xi, \varphi_0(x) \models A(\varphi_{1x}(a_1), \ldots \varphi_{1x}(a_m)) \bigr )
	\]
	
	We consider only two cases. The other cases are straightforward.
	
	\textbf{Case 1.} $ A = \Box B $.
	Let $ \pX, \theta, x \models \Box B $, then 
	$$\exists U \in \tau(x) \forall y \in U (\pX, \theta, y \models B).  $$  
	For each $ k \in \set{1, \ldots, m} $ there exists $ U_k \in \tau(x)$  such that $ \varphi_{1y}(a_k) = \varphi_{1x}(a_k)$ for any  $ y \in U_k $. 
	
	We take $ U' = U \cap \bigcap_{k=1}^{m} U_k $. Set $ U' \in  \tau(x)$ since $ \tau(x) $ is a filter. Note that $ B $ is true at all points in $ U' $ and by the induction hypotheses $ B(a_1, \ldots a_m) $ is true at all points in $ \varphi_0(U')$.
	
	By the definition of p-morphism of n-frames $ \varphi_0(U') $ is a neighbourhood of point $ \varphi_0(x) $.  Since $R(\varphi_0(x))$ is the minimum neighbourhood of point $\varphi_0(x)$ in frame $ \mN(W, R) $ then $R(\varphi_0(x)) \subseteq \varphi_0(U')  $. Hence, $ \pF, \xi, u \models B(\varphi_{1y}(a_1), \ldots \varphi_{1y}(a_m)) $ for any $ u \in R(\varphi_0(x))$, and any $ y \in U'$ such that $ \varphi_0(y) = u $. Since for any $ y\in U' $ $\varphi_{1y}(a_1) = \varphi_{1x}(a_1), \ldots \varphi_{1y}(a_m) = \varphi_{1x}(a_m)$ then 
	$$
	\pF, \xi, \varphi_0(x) \models \Box B (\varphi_{1x}(a_1), \ldots \varphi_{1x}(a_m)).
	$$
	
	Now let $ \pF, \xi, \varphi_0(x) \models \Box B (\varphi_{1x}(a_1), \ldots \varphi_{1x}(a_m))$, hence 
	$$ 
	\pF, \xi, u \models B (\varphi_{1x}(a_1), \ldots \varphi_{1x}(a_m))\hbox{ for all } u \in R(\varphi_0(x)).
	$$ 
	By the definition of p-morphism of n-frames $ U = \varphi_0^{-1}(R(\varphi_0(x))) \in \tau(x)$. 
	For each $ k \in \set{1, \ldots, m} $ there exists $ U_k \in \tau(x)$  such that $ \varphi_{1y}(a_k) = \varphi_{1x}(a_k) $ for any  $ y \in U_k$. We take $ U' = U \cap \bigcap_{k=1}^{m} U_k \in \tau(x)$. Then  $ \pF, \xi, u \models B (\varphi_{1y}(a_1), \ldots \varphi_{1y}(a_m))$ for all $ u \in R(\varphi_0(x))$ and any $ y \in U'$ such that $ \varphi_0(y) = u $.
	By induction hypotheses 
	$$ 
	\forall y \in U' \left (\pX, \theta, y \models B(a_1, \ldots a_m)\right ),\hbox{ hence } \pX, \theta, x \models \Box B (a_1, \ldots a_m).
	$$
	
	\textbf{Case 2.} $ A = \forall t B(t, a_1, \ldots a_m) $. 
	$ \pX, \theta, x \nvDash \forall t B(t, a_1, \ldots a_m) $ then there exists $ a\in D^* $ such that $ \pX, \theta, x \nvDash B(a, a_1, \ldots a_m) $. By induction hypotheses 
	\[ 
	\pF, \xi, \varphi_0(x) \nvDash  B(\varphi_{1x}(a), \varphi_{1x}(a_1), \ldots \varphi_{1x}(a_m)), \]
	hence $ \pF, \xi, \varphi_0(x) \nvDash \forall t B(t, \varphi_{1x}(a_1), \ldots \varphi_{1x}(a_m)) $.
	
	Let $ \pF, \xi, \varphi_0(x) \nvDash \forall t B(t, \varphi_{1x}(a_1), \ldots \varphi_{1x}(a_m)) $ then there exists $ b\in D_{\varphi_0(x)} $ such that $ \pF, \xi, \varphi_0(x) \nvDash B(b, \varphi_{1x}(a_1), \ldots \varphi_{1x}(a_m)) $. Since $ \varphi_{1x} $ is surjective there exists $ a \in D^* $ such that $ \varphi_{1x}(a) = b $, then by induction hypotheses $ \pX, \theta, x \nvDash  B(a, a_1, \ldots a_m) $, hence $ \pX, \theta, x \nvDash \forall t B(t, a_1, \ldots a_m) $.
\end{proof}

The following lemma states that the composition of two $ p $-morphisms is a $ p $-morphism. 
\begin{lem}\label{lem:pmorphism_composition}
	Let $ \pX $ be an n-frame with constant domain, $ \pF_1 $ and $ \pF_2 $ be two Kripke frames with expanding domains. 
	If $ (\varphi_0, \varphi_1):\pX \pmor \pF_1 $ and $ (\psi_0, \psi_1):\pF_1 \pmor \pF_2 $ then 
	$$ (\psi_0 \circ \varphi_0, \eta):\pX \pmor \pF_2 ,$$
	where $ \eta = (\eta_x)_{x\in X} $,  $ \eta_x(d) = \psi_{1\varphi_0(x)}(\varphi_{1x}(d))$ for $ d\in D^* $.
\end{lem}
The proof is straightforward.

\section{Main construction}

%To simplify the construction we assume that $ N = 1$, i.e.{} we have only one box and one relation in Kripke frames. Later we explain how to adapt this construction to multiple modalities.   

The following construction was introduced in \cite{kudinov_aiml12,kudinov2014neighbourhood}. 

For a Kripke frame $ F $ Starting from a Kripke frame we will construct a neighbourhood frame no point has the smallest neighbourhood. In topology a topological space is dense if it has no isolated points.

\begin{definition}
	Let $ \Sigma $ be a non-empty set (\emph{alphabet}). A finite sequence of elements from $ \Sigma $ is a \emph{word}; by $\eps$ we denote the empty word.     
	Let $\Sigma^*$ be the set of all words. We will write words without brackets or commas, e.g.{} $ a_1 a_2 \ldots a_n \in \Sigma^* $. The length of a word is the number of elements in it:
	\[ 
	\len(a_1 a_2 \ldots a_n) = n, \quad \len(\eps)=0.
	\]
	
	The \emph{concatenation} of words is defined as follows:
	\[ 
	a_1 a_2 \ldots a_n \cdot b_1 b_2 \ldots b_m = a_1 a_2 \ldots a_n  b_1 b_2 \ldots b_m
	\] 
\end{definition}

\begin{definition}
	For a frame $F=(W, R)$ with root $a_0$ %R_1, \ldots, R_N)$ 
	we define \emph{a (rooted) path with stops} as a word in alphabet $W \cup \set{0}$: $a_1 \ldots a_n$, so that $a_i \in W$ or $a_i = 0$,  and after dropping zeros, each point is related to the next one by relation $R$, and the first one is a successor of the root. 
	The empty word $\eps$ is allowed.
	Any path with stops is a word of the following type
	\begin{align*}
		&0^{i_0} b_1 0^{i_1} \ldots 0^{i_{m-1}} b_m 0^{i_{m}},\ \hbox{where } b_j \in W,\ i_j \ge 0,\ 0^i = \underbrace{00\ldots 0}_{\hbox{\small $i$ times}};\\
		\hbox{and }  &f_0(0^{i_0} b_1 0^{i_1} \ldots b_m 0^{i_{m}}) = a_0 R b_1 R \ldots R b_m \in W^{\sharp},\\
		&f_0(\eps) = f_0(0^{i_0}) = a_0 \in W^{\sharp},
	\end{align*}
	where $ W^{\sharp} $ is the set of all rooted paths in $ F $ (Definition \ref{def:unraveling})
	
	Let us consider some examples. 
	\begin{itemize}
		\item $f_0(\eps) = a_0$.
		\item If the root $ a_0 $ is reflexive, then $a_0$, $ a_0 000a_0 $, $ a_0a_0a_0 $ are examples of paths with stops, and $f_0(a_0) = a_0 R a_0$. 
		\item In frame $F = (W,R)$, where $W=\set{a_0,b}$, $R=\set{(a_0,b)}$ the following words are paths with stops: $\eps$, $000$, $00b$, $00b00$. In general any path with stops in $ F $ equals to $ 0^k b 0^m $ or to $ 0^k $ for some $ k,m \ge 0 $.
		\item In frame $F' = (W,R')$, where $W=\set{a_0,b}$, $R'=\set{(a_0,b), (b, a_0)}$ the following words are paths with stops: $\eps$, $000$, $00b$, $00b00a_0$, $ b a_0 b a_0 $. But $ a_0 $ is not a path with stops. 
		
	\end{itemize}

	%the relations $R_1,\ \ldots, R_N$.
	We also consider infinite paths with stops that end with infinitely many zeros. We call these sequences \emph{pseudo-infinite paths (with stops)}. 
	A pseudo-infinite path can be presented uniquely in the following way:
	\[
	\alpha = 0^{i_1} b_1 0^{i_2} \ldots 0^{i_m} b_m 0^{\omega},\\ \hbox{where } b_j \in W,\ i_j \ge 0.
	\]
	Let $W_\omega$ be the set of all pseudo-infinite paths in $W$.
\end{definition}

In the following we define function $f_0:W_\omega \to W^\sharp$. For a pseudo-infinite path $\alpha = a_1 \ldots a_n \ldots $ we put
\begin{align*}
	\st(\alpha) &= \min \setdef[m]{\forall k > m (a_k = 0)};\\
	\alpha|_k &= a_1 \ldots a_k;\ \ \alpha|_0 = \eps;	\\ 
	f_0(\alpha) &= f_0(\alpha|_{\st(\alpha)}), i.e.,\ f_0(0^{i_1} b_1 0^{i_2} \ldots 0^{i_m} b_m 0^{\omega}) = a_0 R b_1 R \ldots R b_m.
\end{align*}

In order to introduce a neighbourhood function on $W_{\omega}$, for $ k\ge 0 $, we define
\[ 
U_k(\alpha) = \Setdef[\beta \in W_\omega]{\alpha |_m = \beta |_m \ \&\ f_0(\alpha) R^\sharp f_0(\beta), \ m = \max(k, \st(\alpha))}.
\]

\begin{lem}\label{lem:U_m-base}
	$U_k (\alpha) \subseteq U_m (\alpha)$ whenever $k \ge m$. %for any $i \in \set{1,2}$. %need to prove in a lemma 
\end{lem}
\begin{proof}
	Let $\beta \in U_k (\alpha)$. Since $\alpha |_k = \beta |_k$ and $k \ge m$,  $\alpha |_m = \beta |_m$. Hence, $\beta \in U_m (\alpha)$.
\end{proof}

\begin{definition}\label{def:n-frame_from_fframe}
	Due to Lemma \ref{lem:U_m-base}, set $\set{U_n(\alpha)}_{n\in \Natr}$ form a filter base. So we can define
	\begin{align*}
		\tau (\alpha) &- \hbox{the filter with the base } \setdef[U_n(\alpha)]{n \in \Natr};\\
		\Nf(F) &= (W_\omega,\tau) - \hbox{is \emph{a dense n-frame based on} $F$.} 
	\end{align*}
\end{definition}

In frame $\mathcal{N_\omega}(F)$ no point has the minimal neighbourhood unlike $\mN(F)$. Indeed,
\begin{equation}  \label{eq:density}
\bigcap\limits_{n\in \Natr} U_n(\alpha) = \varnothing \not\in \tau(\alpha). 
\end{equation}

To prove (\ref{eq:density}) let us take $ \beta \in U_m(\alpha) $, then $ \beta = \alpha|_k b' 0^\omega $ for some $ k\ge m $ and $ b' \in W $. But then $ \beta \notin U_{k+1}(\alpha) $. 

\begin{lem}
	Equality  $ f_0(U_k(\alpha)) = R^{\sharp}(f_0(\alpha)) $ holds for any $ \alpha $ and $ k $. 
\end{lem}

\begin{proof}
	For $ k<st(\alpha) $ $ U_k(\alpha) = U_{st(\alpha)}(\alpha) $, hence we can assume that $ k
	\ge st(\alpha) $.
	
	Let $ x \in  f_0(U_k(\alpha)) $. Then 
	$$ 
	\exists \beta (\alpha |_k = \beta |_k \ \&\ f_0(\alpha) R^\sharp f_0(\beta) = x),$$ 
	so $ x\in R^{\sharp}(f_0(\alpha)) $.
	
	Let $ x \in R^{\sharp}(f_0(\alpha)) $, then 
	$x = f_0(\alpha) R b'$ for some $ b' $.
	Consider
	\[ 
	\beta = \alpha 0^{k-st(\alpha)} b'.
	\]
	It is easy to show that $ \beta \in  U_k(\alpha)$ and $ f_0(\beta) = x $.
\end{proof}

\begin{lem}\label{lem:bmorphism_wframe2nframe}
	Let $F = (W, R)$ be a Kripke frame with root $a_0$, then
	\[
	f_0 : \Nf(F) \pmor \mN(F^\sharp).
	\]
\end{lem}

\begin{proof}
	Since for any $b \in W$ there is a path $a_0 R a_1 R\ldots R a_{n-1} R b$, hence for a pseudo-infinite path $\alpha =  a_1\ldots b 0^\omega \in X$, $f(\alpha) = b$. So $f_0$ is surjective.
	
	Let us prove the zig property. Assume that $\alpha\in W_\omega$ and $U \in \tau(\alpha)$. We need to prove that $R^\sharp(f_0(\alpha)) \subseteq f_0(U)$. 
	There exists $m$ such that $U_m(\alpha) \subseteq U$, and since $f_0(U_m(\alpha)) = R^\sharp(f_0(\alpha))$, we have
	\[
	R^\sharp(f_0(\alpha)) = f_0(U_m(\alpha)) \subseteq f_0(U).
	\]
	
	Let us prove the zag property. Assume that $\alpha\in W_\omega$ and $V$ is a neighbourhood of $f_0(\alpha)$, i.e.{} $R^\sharp(f_0(\alpha)) \subseteq V$. We need to prove that there exists $U \in \tau(\alpha)$, such that $f_0(U) \subseteq V$. As $U$ we take $U_m(\alpha)$ for some $m \ge st(\alpha)$, then
	\[
	f_0(U_m(\alpha))=R^\sharp(f_0(\alpha)) \subseteq V. 
	\]
\end{proof}

\begin{cor}\label{cor:logicNf}
	For any frame $F$, $\Log(\Nf(F)) \subseteq \Log(F)$.
\end{cor}

\begin{proof}	
	It follows from Lemmas \ref{lem:n-frame_from_kframe},  \ref{lem:bmorphism_wframe2nframe}, \ref{lem:unrev_pmorph}, and Corollary \ref{cor:pmorph} that
	\[
	\Log(\Nf(F)) \subseteq \Log(\mN(F^\sharp)) = \Log(F^\sharp) \subseteq \Log(F). 
	\]
\end{proof}

Let us remark that it is possible that $ \Log(\Nf(F)) \ne \Log(F) $. For example, consider the natural numbers with the ``next" relation. It is convenient here to regard a number as a word in a one-letter alphabet: 
\[ 
G = (\set{1}^* , S),\ 1^n S 1^m \iff m=n+1.
\] 

Obviously $ G \models \romb p \to \Box p $.

Since in $ G $ every point, except for the root, has only one predecessor, we can identify a point and a path from the root to this point, i.e.,\ $ G^\sharp = G $. Therefore, points in $\Nf(G)$ can be presented as infinite sequences of 0 and 1 with only zeros at the end. 
\begin{prop}
	$ \Nf(G) \nvDash \romb p \to \Box p $
\end{prop}
\begin{proof}
	Consider valuation $ V(p) = \setdef[0^{2n}10^{\omega}]{n \in \Natr} $. In every neighbourhood of point $ 0^\omega $ there are points where $ p $ is true and there are points where $ p $ is false. 
	For example, if $ k $ is even, then in $U_k(0^{\omega})$ there is a point $ 0^k 1 0^{\omega} $ where $ p $ is true, and a point $ 0^{k+1} 1 0^{\omega} $ where $ p $ is false.
	Hence,
	\[ 
	\Nf(G) \models \romb p \land \romb\lnot p. 
	\]  
\end{proof}

%This formula does not preserved under the $ \Nf $ operation and probably any formula that restricts branching does not preserved under the $ \Nf $ operation.

\begin{definition}
	Let $F_1 = (W_1, R_1)$ and $F_2= (W_2, R_2)$ be two Kripke frames with roots $x_0$ and $y_0$ respectively. Let $ \Sigma=W_1 \cup W_2 $ then we define functions $ p_1: \Sigma^* \to W_1^* $, $ p_2: \Sigma^* \to W_2^* $ and $ \pi: \Sigma^*\setminus \set{\eps} \to \Sigma$ by induction
	\begin{equation*}
		\begin{array}{rll}
			p_1 (\eps) &= x_0, \\
			p_2 (\eps) &= y_0, \\
			\pi (u) &= u,& \hbox{ for } u \in \Sigma,\\
			p_1 (\vec a u) &= p_1(\vec a)\cdot u,&\hbox{ for } \vec a \in \Sigma^*,\ u \in W_1, \\
			p_1 (\vec a u) &= p_2(\vec a),&\hbox{ for } \vec a \in \Sigma^*,\ u \in W_2,\\
			p_2 (\vec a u) &= p_1(\vec a),&\hbox{ for } \vec a \in \Sigma^*,\ u \in W_1, \\
			p_2 (\vec a u) &= p_2(\vec a)\cdot u,&\hbox{ for } \vec a \in \Sigma^*,\ u \in W_2,\\
			\pi(\vec a u) &= u, &\hbox{ for } \vec a \in \Sigma^*,\ u \in \Sigma.
		\end{array}
	\end{equation*}
\end{definition}
Since $ F_1$ and $ F_2$ are Kripke frames with one relation we can assume that paths in them do not contain relations and it will not lead to misunderstanding:
\begin{align*}
	W_1^\sharp &= \setdef[x_0 x_1\ldots x_n]{x_0 R_1 x_1 R_1\ldots R_1 x_n \hbox{ --- a path in the usual sense}}\\    
	W_2^\sharp &= \setdef[y_0 y_1 \ldots y_n]{y_0 R_2 y_1 R_2 \ldots R_2 y_n \hbox{ --- a path in the usual sense}}\\
\end{align*}

We define the \emph{entanglement} of $ F_1 $ and $ F_2 $ as follows
\begin{align*}
	F_1 \mathop{\taurus}  F_2 &= \setdef[\vec{x} \in \Sigma^*]{ p_1(\vec{x}) \in W_1^{\sharp} \hbox{ and } p_2(\vec{x}) \in W_2^{\sharp}},\\
	\vec a \mathop{\taurus} F_2 &= \setdef[\vec{x} \in F_1 \taurus  F_2]{ p_1(\vec{x}) = \vec a }, \hbox{ for } \vec a \in W^\sharp_1.\\
\end{align*}

A universal quantifier behave a lot like an $ \logic{S5} $ modality. So we use entanglement with an $\logic{S5}$-frame to construct a predicate frame.

Let $ F = (W, R) $ be a rooted propositional Kripke frame and $ G = (\Real_{\emptyset}, \Real_{\emptyset}\times \Real_{\emptyset}) $ be the continuum $ \logic{S5} $-frame, here $ \Real_{\emptyset} = \Real \setminus \set{0} $ and $ \Real $ is the set of all real numbers. As the underling frame we take $ F^\sharp $ and for a path $ \vec a \in W^\sharp$ we take the following set $ D'_{\vec a} =  \vec a \mathop{\taurus} G$. But this cannot be our family of domains since for any $ \vec a \ne \vec b $ we have $ D'_{\vec a} \cap D'_{\vec b} = \varnothing$. Let us define the following equivalence relation on $ F \mathop{\taurus} G $. For $ \vec x, \vec y \in  F \mathop{\taurus} G$ 
\[ 
\vec x \sim \vec y \iff \exists \vec t \in  F \mathop{\taurus} G \;\exists \vec c, \vec d \in W^* (\vec x  = \vec t \cdot \vec c \hbox{ and } \vec y = \vec t \cdot \vec d ).
\]
It is easy to check that $ \sim $ is reflexive and symmetric.
Let us check the transitivity.
\begin{align*}
	\vec x \sim \vec y \ \&\ \vec y \sim \vec z &\Rightarrow \\
	\exists \vec t_1 \in  F \mathop{\taurus} G \;\exists \vec c_1, \vec d_1 \in W^* &(\vec x  = \vec t_1 \cdot \vec c_1 \ \& \ \vec y = \vec t_1 \cdot \vec d_1 ) \ \&\\
	\exists \vec t_2 \in  F \mathop{\taurus} G \;\exists \vec c_2, \vec d_2 \in W^* &(\vec y  = \vec t_2 \cdot \vec c_2 \ \& \ \vec z = \vec t_2 \cdot \vec d_2 )
\end{align*}
Since $ \vec y = \vec t_1 \cdot \vec d_1 = \vec t_2 \cdot \vec c_2 $ then there exist $ \vec e \in W^* $ such that $ \vec y = \vec t_1 \cdot \vec e \cdot \vec c_2$ or $ \vec y = \vec t_2 \cdot \vec e \cdot \vec d_1$. 

If $ \vec y = \vec t_1 \cdot \vec e \cdot \vec c_2$ then $ \vec z = \vec t_1 \cdot \vec e \cdot \vec d_2 $ and $ \vec x \sim \vec z $.

If $ \vec y = \vec t_2 \cdot \vec e \cdot \vec d_1$ then $ \vec x = \vec t_2 \cdot \vec e \cdot \vec c_1 $ and $ \vec x \sim \vec z $.

Therefor $ \sim $ is an equivalence relation.

Let $ [\vec x] $ be the equivalence class of $ \vec x $. We define
\begin{align}
	D^+ &= F \mathop{\taurus} G / {\sim}  = \setdef[{[\vec x]}]{\vec x \in F \mathop{\taurus} G};\\
	D^\sharp_{\vec a} &= \setdef[{[\vec x]}]{\vec x \in D'_{\vec a}}. \label{eq:Dsharpvec}
\end{align}

\begin{lem}
	For $ \vec a $ and $ \vec b $ from $ F^\sharp $ if $ \vec a R^{\sharp} \vec b $ then $ D^\sharp_{\vec a} \subset D^\sharp_{\vec b} $.
\end{lem}

\begin{proof}
	Since $ \vec a R^{\sharp} \vec b $ then there exists $ \vec c \in F^\sharp $ such that $ \vec a \cdot \vec c = \vec b $.
	So
	\begin{align*} 
		[\vec x] \in D^\sharp_{\vec a} \;\Leftrightarrow&\; \exists \vec y \in D'_{\vec a} (\vec x \sim \vec y)
		\Rightarrow \vec y \cdot \vec c \in D'_{\vec b}
	\end{align*}
	But $ \vec y \cdot \vec c \sim \vec y \sim \vec x$, hence $ [\vec x ]\in D^\sharp_{\vec b}$.
\end{proof}

Let $ \pF = (F, D) $ be a predicate Kripke frame, $ F = (W, R) $ be a Kripke frame, and cardinality of all $ D_u $ ($ u\in F $) is no greater than continuum. 

We take $ \pF^\sharp = (F^\sharp, D^\sharp) $, where $ D^\sharp = (D^\sharp_{\vec a})_{\vec a \in W^\sharp} $ and $ D^\sharp_{\vec a} $ is from (\ref{eq:Dsharpvec}).

Note that $ \pi $ restricted to $ F^\sharp $ is a p-morphism from $ F^\sharp $ to $ F $ (c.f.{} Lemma \ref{lem:unrev_pmorph}). 

Let us define $ \psi = (\psi_{\vec a})_{\vec a \in W^\sharp} $ by induction:
\begin{enumerate}
	\item Let $ \psi_{\eps} $ be an arbitrary surjective map from $ D^\sharp_{\eps} $ to $ D_{a_0} $, where $ a_0 $ is the root of $ F $. Such map exists since cardinality of $ D^\sharp_{\eps} $  is continuum.
	\item Assume that $ \psi_{\vec a} $ is already defined and $ \vec b = \vec a u $. 
	Note that cardinality of set $ D^\sharp_{\vec b} \setminus D^\sharp_{\vec a} $ is continuum, so there is a surjective map $ \eta $ from $ D^\sharp_{\vec b} \setminus D^\sharp_{\vec a} $ to $ \left (D_{\pi(\vec b)} \setminus D_{\pi(\vec a)}\right ) \cup \set{[\eps]}$. 
	We need to add $ \set{[\eps]} $, since it is possible that $ D_{\pi(\vec b)} = D_{\pi(\vec a)} $. We put 
	\[
	\psi_{\vec b}([\vec x]) = \begin{cases}
		\psi_{\vec a}([\vec x]),\hbox{ if } [\vec x] \in {D^\sharp_{\vec a}};\\
		\eta([\vec x]), \hbox{ otherwise.}
	\end{cases}
	\]
\end{enumerate}

\begin{lem}\label{lem:pmorphism:unrevelpF_to_pF}
	If frame $ F $ is a tree, i.e.{} there is a unique path from the root to any other point, then $ (\pi, \psi):\pF^\sharp \pmor \pF $.
\end{lem}

To prove this we should check all items of the definition of p-morphism, and it is an easy exercise.  

Now we define our predicate neighbourhood frame.
\begin{align*}
	D^* &= \Nf (G);\\
	\X  &= \Nf (F);\\
	\pX &= (\X, D^*).
\end{align*}

By Lemma \ref{lem:bmorphism_wframe2nframe} $ f_0: \X \pmor \mN(F^\sharp) $.

Let us define functions $ \xi_\alpha:D^* \to D^\sharp_{f_0(\alpha)} $ for each $ \alpha \in \X $. Function $ \xi_\alpha $ applied to $ \gamma \in D^* $ replaces zeros in $ \alpha $ with letters from $ \gamma $. To define it properly we first define function $ h:W_\omega \times \Real_\omega \to (W \cup \Real_{\emptyset})^* $ by induction on $ \st(\alpha),$ here $ \alpha \in W_\omega $ and $\gamma \in \Real_\omega$:
\begin{description}
	\item[Base.] $ \st(\alpha) = 0 $: $ \alpha = 0^\omega $, $ h(0^\omega, \gamma) = f_0(\gamma)$;
	\item[Step.] Assume that $ \alpha = a_1 \beta$ and $ \gamma = c_1\delta $. Then
	\[ 
	h(a_1 \beta, c_1\delta) = \begin{cases}
		f_0(a_1 \cdot h(\beta, \delta)),\hbox{ if } c_1 = 0;\\
		c_1 \cdot h(\beta, a_1\delta), \hbox{otherwise}.
	\end{cases}
	\]
\end{description}

For example 
\begin{align*}
	h(a0b00c0^\omega, &\;10340^\omega)= 1\cdot h(a0b00c0^\omega, 0340^\omega) = 1a\cdot h(0b00c0^\omega, 0340^\omega) \\
	&= 1a\cdot h(b00c0^\omega, 340^\omega)=1a3\cdot h(b00c0^\omega, 40^\omega)=\ldots=1a34bc.
\end{align*} 

Now we can define $ \xi_\alpha $:
\[ 
\xi_\alpha(\gamma) = [h(\alpha, \gamma)].
\]

\begin{lem}\label{lem:X-F_pred_p-morphism}
	Pair of functions $(f_0, \xi)  $, where $ \xi = (\xi_\alpha)_{\alpha\in \X} $ is a p-morphism from $ \pX $ to $ \pF^\sharp $.
\end{lem}

\begin{proof}
	In the following we check all the items of Definition \ref{def:pmorph4prednframe}.
	
	Item (\ref{it:1:def:pmorph4prednframe}) follows from Lemma \ref{lem:bmorphism_wframe2nframe}.
	
	Item (\ref{it:2:def:pmorph4prednframe}). We fix $ \alpha\in \X $ and take $ [\vec x] \in D^\sharp_{f_0(\alpha)} = D'_{f_0(\alpha)} / {\sim}$. Consider $ \vec y \in D'_{f_0(\alpha)} =  f_0(\alpha) \taurus G$, such that $ \vec y \in [\vec x] $. We should insert in the sequence $ \alpha $ numbers from $ \Real_\emptyset $ so that $ f_0(\alpha') = \vec y $ where $ \alpha' $ is the resulting sequence. After that we replace all symbols from $ W $ in $ \alpha $ by zeros and get $ \gamma \in \Real_\omega $. And $ \xi_{\alpha}(\gamma) = [\vec y] = [\vec x]$.

To show that it is always possible, we define function $ t $ by induction on $ st(\alpha) + \len(\vec y)$:
\begin{align*}
	t(0^\omega, \vec{y}) &= \vec{y}\cdot 0^\omega\\
	t(a_1\beta, c_1 \vec{z}) &= \begin{cases}
		a_1\cdot t(\beta, c_1 \vec{z}),\hbox{ if } c_1 = 0 \ \&\ a_1 ;\\
		a_1\cdot t(\beta, \vec{z}),\hbox{ if } c_1 \ \& \ a_1 \in W;\\
		c_1 \cdot t(a_1 \beta, \vec{z}), \hbox{ if } c_1 \in \Real_\emptyset.\\
	\end{cases}
\end{align*}

	We illustrate this by an example:
	\begin{align*}
		\alpha &= ab00c0^\omega &		f_0(\alpha) &= abc \\
		\vec{y} &=a12bc3 & \alpha' &= t(\alpha, \vec y) = a12b00c30^\omega
	\end{align*} 
Indeed,		
	\begin{align*}
		&t(ab00c0^\omega, a12bc3) = a\cdot t(b00c0^\omega, 12bc3) = a1\cdot t(b00c0^\omega, 2bc3) = \\
		=\ &a12\cdot t(b00c0^\omega, bc3) = a12b\cdot t(00c0^\omega, c3) = a12b0\cdot t(0c0^\omega, c3) = \\
		=\ &a12b00\cdot t(c0^\omega, c3) = a12b00c\cdot t(0^\omega, 3) = a12b00c3\\
		\gamma =\ &01200030^\omega \qquad \xi_\alpha(\gamma) = [h(\alpha, \gamma)] = [a12bc3]
	\end{align*} 

Let us show that $ f_0(t(\alpha, \vec y)) = \vec y $. For the base of induction we have  
$$ f_0(t(0^\omega, \vec{y})) = f_0(\vec{y}\cdot 0^\omega) = \vec{y}.$$ 
Let $ \alpha = a_1\beta $ and $ \vec{y} = c_1\vec{z} $
\begin{itemize}
	\item [(case 1)] $ c_1 = 0 $ and $ a_1 = 0 $, $ f_0(t(a_1\beta, c_1 \vec{z})) = f_0(0\cdot t(\beta, \vec{y})) = f_0 (t(\beta, \vec{y})) = \vec{y}$;
	\item [(case 2)] $ c_1 \in W $ and $ a_1 \in W $ (note, that in this case $ a_1 = c_1 $ because otherwise $\vec{y} \notin f_0(\alpha) \taurus G $), 
	$$  
	f_0(t(a_1\beta, c_1 \vec{z})) =  a_1\cdot f_0(t(\beta, \vec{z})) = a_1 \cdot f_0(t(\beta, \vec{z})) = a_1 \cdot \vec{z} = \vec{y};
	$$
	\item [(case 3)] $c_1 \in \Real_\emptyset  $, $ f_0(t(a_1\beta, c_1 \vec{z})) = f_0(c_1 \cdot t(a_1 \beta, \vec{z})) = c_1 \cdot f_0(t(a_1 \beta, \vec{z})) = c_1 \vec{z} = \vec{y}$.
\end{itemize}

	Item (\ref{it:3:def:pmorph4prednframe}). 	
 Consider $ \gamma\in \Real_\omega $, $ \alpha \in \X $ and $ m = \st(\gamma) + \st(\alpha)$. If $ \beta \in U_m(\alpha) $ then
	\[ 
	\xi_{\beta}(\gamma) = [h(\beta, \gamma)] = [h(\rest{\alpha}_{st(\alpha)}0^kb0^\omega, \gamma)] = [h(\alpha, \gamma)\cdot b] = [h(\alpha, \gamma)] = \xi_{\alpha}(\gamma).
	\]	
	
	In our example$ m = st(\gamma) + st(\alpha) = 11$ and $ \beta \in  U_m(\alpha)$ then $ \beta = ab0c0000000 \cdot \vec z \cdot 0^\omega $ for some $ \vec z \in (W \cup \set{0})^*$. Then
	\begin{align*}
		\xi_{\alpha}(\gamma) &= [a12bc3]\\
		\xi_{\beta}(\gamma) &= [a12bc3 \cdot f_0(\vec z)] = [a12bc3]
	\end{align*} 
	
	This finishes the proof. 
\end{proof}

\section{Completeness results }
This already gives us the following theorem 
\begin{thm}[\cite{arlo2006first}]
	Logic $ \logic{QK} $ is complete with respect to neighbourhood frames with constant domain.
\end{thm}

Indeed, if $ A \notin \logic{QK} $ then by Theorem \ref{thm:PTC_ncompleteness} it falsifiable in a Kripke frame $ \pF $ with expanding domains. By Lemmas \ref{lem:pmorphism:unrevelpF_to_pF}, \ref{lem:X-F_pred_p-morphism}, \ref{lem:pmorphism_composition} and \ref{lem:pmorphismlemma4nF2kF} $ \Log(\pX) \subseteq\Log(\pF)  $. So $ A \notin \Log(\pX) $.

\begin{definition}
	Let $\Gamma$ be a set of universal strict Horn sentences, $F = (W,R)$ be a rooted frame, $\alpha \in W_{\omega}$ and $f_0:W_{\omega} \to W^{\sharp}$ be the ``zero-dropping'' function. Then we define 
	\begin{align*}
		U_k^{\Gamma}(\alpha) &= \Setdef[\beta \in W_\omega]{\alpha |_m = \beta |_m \ \&\ f_0(\alpha) (R^\sharp)^{\Gamma} f_0(\beta), \ m = \max(k, \st(\alpha))};\\ 
		\tau^{\Gamma}(\alpha) &= \setdef[V]{\exists k \left ( U_k^{\Gamma}(\alpha) \subseteq V \right )};\\
		\mN^{\Gamma}_{\omega}(F) &= (W_{\omega}, \tau^{\Gamma}).
	\end{align*}
\end{definition}

\begin{lem}
	Let $\Lo$ be a one-way PTC-logic, $\Gamma$ be the corresponding set of Horn sentences, and $F\models \Lo$. If $ \Box p \to \Box^n p \in \Lo $, then  
	\[  
	\mN^{\Gamma}_{\omega}(F) \models \Box p \to \Box^n p.
	\] 
\end{lem}

\begin{proof}  
	Let $ M = (\mN^{\Gamma}_{\omega}(F), V)$ be a neighbourhood model.
	We assume that $M, \alpha \not\models \Box^n p$, and then we prove that
	$M, \alpha \not\models \Box p$, i.e.,
	\[ 
	\forall m \exists \beta \in U^{\Gamma}_m(\alpha)(\beta \not\models p).
	\]
	Let us fix $m$. Then
	\begin{align*}
		&\exists  \alpha_1\in U^{\Gamma}_{m}(\alpha) \left( 
		\alpha_1 \not\models \Box^{n-1} p \right)  \\
		\Rightarrow &\exists \alpha_2\in U^{\Gamma}_{m}(\alpha_1) \left( 
		\alpha_2 \not\models \Box^{n-2} p \right) \\
		\qquad&\vdots\qquad\vdots\qquad\vdots\qquad\vdots\qquad\vdots\qquad\vdots \\
		\Rightarrow &\exists \alpha_n\in U^{\Gamma}_m(\alpha_{n-1}) \left( 
		\alpha_n \not\models p \right). \\
	\end{align*}
	%Remember that $ f_0: W_{\omega} \to W^{\sharp} $ is the forgetting-zeros-function. 
	By the definition of $ U^{\Gamma}_m(\alpha)$
	\[ 
	f_0(\alpha) (R^\sharp)^{\Gamma} f_0(\alpha_1) (R^\sharp)^{\Gamma} \ldots, (R^\sharp)^{\Gamma} f_0(\alpha_n)
	\]
	and
	\[ 
	\rest{\alpha}_m = \rest{\alpha_1}_m = \ldots,  =\rest{\alpha_n}_m.
	\]
	Since $ \left (W^{\sharp}, (R^\sharp)^{\Gamma}\right ) \models \Box p \to \Box^n p$, then
	\[ 
	f_0(\alpha) (R^\sharp)^{\Gamma} f_0(\alpha_n).
	\]
	It follows that $ \alpha_n \in U^{\Gamma}_{m}(\alpha) $.
\end{proof}

\begin{lem}
	Let $\Lo$ be a one-way PTC-logic, $\Gamma$ be the corresponding set of Horn sentence, and $F\models \Lo$. Then
	\[ 
	f_0: \mN^{\Gamma}_{\omega}(F) \pmor \mN(F^{\sharp\Gamma}).
	\]    
\end{lem}
\begin{proof}
	The surjectivity was established in Lemma \ref{lem:bmorphism_wframe2nframe}.
	
	Assume that $\alpha\in W_\omega$ and $U \in \tau^\Gamma(\alpha)$. We need to prove that $R^{\sharp\Gamma}(f_0(\alpha)) \subseteq f_0(U)$. 
	There exists $m$ such that $U^{\Gamma}_m(\alpha) \subseteq U$. 
	
	It is easy to check by the definition that $f_0(U^\Gamma_m(\alpha)) = R^{\sharp\Gamma}(f(\alpha))$. Then
	\[
	R^{\sharp\Gamma}(f_0(\alpha)) = f_0(U^\Gamma_m(\alpha)) \subseteq f_0(U).
	\]
	
	Assume that $\alpha\in W_\omega$ and $U'$ is a neighbourhood of $f_0(\alpha)$, i.e.\ $R^{\sharp\Gamma}(f_0(\alpha)) \subseteq U'$. We need to prove that there exists $U \in \tau^\Gamma(\alpha)$ such that $f(U) \subseteq U'$. Let us take $U = U^{\Gamma}_m(\alpha)$ for some $m \ge st(\alpha)$, then
	\[
	f_0(U^\Gamma_m(\alpha))=R^{\sharp\Gamma}(f_0(\alpha)) \subseteq U'. 
	\]
\end{proof}

Let $ F $ be a rooted $ \logic{L} $-frame. Let us define 
\begin{align*}
	\X^\Gamma &= \Nf^\Gamma(F), \\ 
	\pX^\Gamma &= (\X^\Gamma, D^*),\hbox{ where } D^* = \Nf(\Real_{\emptyset}, \Real^2_{\emptyset}),\\
	\pF^{\sharp\Gamma}&=(F^{\sharp\Gamma}, D^\sharp).
\end{align*}

\begin{lem}
	$ \pF^{\sharp\Gamma} $ is a predicate Kripke frame with expanding domains.
\end{lem}

We need to check that if $ \vec a R^\sharp \vec b $ then $ D^\sharp_{\vec a} \subseteq  D^\sharp_{\vec b} $. This is true since relation $ \subseteq $ is transitive and $ \Gamma $-closure is included in the transitive-reflexive closure.

The underlying sets of $ \X $ and $ \X^\Gamma $ are the same, so function $ \xi $ defined in the previous section will be well-defined here as well.
\begin{lem}
	Pair $ (f_0, \xi) $ is a p-morphism from $ \pX^\Gamma $ to $ \pF^{\sharp\Gamma} $.
\end{lem}

The proof is the same as in Lemma \ref{lem:X-F_pred_p-morphism}. For the last condition we need to check that for any $ \alpha, \beta \in \X^\Gamma $ such that $ \beta \in U^\Gamma_m(\alpha) $ and any $ \gamma \in D^* $ 
\[ 
\xi_{\beta}(\gamma) = \xi_{\alpha}(\gamma).
\]
This is an easy exercise.

\begin{thm}\label{thm:main_PTClogic-nframecompleteness}
	Let $ \logic{L} $ be a one-way PTC-logic with one modality then predicate modal logic $ \logic{QL} $ is complete with respect to predicate neighbourhood frames with constant domain. 
\end{thm}

\section{Multiple modalities}

The construction should be working for multiple modalities if as the alphabet for sequences we take not just elements of $ W $ ($ W_1 $ or $ W_2 $) but pairs $ (R_i, w) $, where $ i\in \set{1, \ldots, N} $ and $ w\in W $, $ N $ is the number of modalities.  All the definitions should be changed accordingly.

The author strongly believes that the following hypothesis can be proven by a straightforward adaptation of the methods from this paper:

\textbf{Hypothesis.} \textit{ Let $ \logic{L} $ be a one-way PTC-logic with arbitrary many modalities then predicate modal logic $ \logic{QL} $ is complete with respect to predicate neighbourhood frames with constant domain.}

\section{Topological semantics}

It is well known that the topological semantics is a particular case of the neighbourhood semantics and an $ \logic{S4} $-neighbourhood frame is basically a topological space with closure operator interpreting the $ \romb $ modality. 

This observation gives us the following previously proven theorem as a corollary of Theorem \ref{thm:main_PTClogic-nframecompleteness} 

\begin{thm}[\cite{rasiowa1963metamathematics}]
	Logic $ \logic{QS4} $ is complete with respect to topological spaces.
\end{thm}

Using construction from \cite{kudinov_aiml12} and Theorem \ref{thm:main_PTClogic-nframecompleteness} we can prove

\begin{thm}[\cite{kremer2014quantified}]
	Logic $ \logic{QS4} $ is complete with respect to the set of rational numbers $ \QQ $.
\end{thm}

As was explained in \cite{kudinov_aiml12} $ \logic{K4^-} = \logic{K} + \Box p \land p \to \Box\Box p $ neighbourhood frame is basically a topological space with derivational operator interpreting the $ \romb $ modality. 
%To be more precise we can define an n-frame based on a topological space such that a set $ U $ is a neighbourhood of point $ x $ if there is an open $ V $ with $ x\in V $ and $ V \setminus \set{x} \subseteq U $. Such n-frame is a \logic{K4^-}-n-frame and truth of formulas preserved by this construction. On the other hand if we have a \logic{K4^-}n-frame $ (X, \tau) $  such that $\forall x\in X\exists U \in \tau(x) (x \notin U) $ then we can construct a topological space $ (X, T) $ such that for any valuation $ V $ for any formula $ A $ $ \forall x\in X (X,\tau, V, x \models A \iff X, T, V, x \models A)$.
From this and the results of the current paper it follows  

\begin{thm}
	Logics $ \logic{QK4} $ and $ \logic{QD4} $ are complete with respect to $ T_d $ and dence-in-itself $ T_d $ topological spaces respectively.
\end{thm}

\begin{thm}
	Logic $ \logic{QD4} $ is complete with respect to the set of rational numbers $ \QQ $ (with derivational modality).
\end{thm}

\bibliographystyle{plain}
\bibliography{quantmodal_arxiv}

\begin{thebibliography}{10}

\bibitem{arlo2006first}
H.~Arl{\'o}-Costa and E.~Pacuit.
\newblock First-order classical modal logic.
\newblock {\em Studia Logica}, 84(2):171--210, 2006.

\bibitem{blackburn_modal_2002}
P.~Blackburn, M.~de~Rijke, and Y.~Venema.
\newblock {\em Modal Logic}.
\newblock Cambridge University Press, 2002.

\bibitem{Chellas1980}
B.~Chellas.
\newblock {\em Modal Logic: An Introduction.}
\newblock Cambridge University Press, Cambridge, 1980.

\bibitem{GSh_Product_1998}
D.~Gabbay and V.~Shehtman.
\newblock Products of modal logics. {P}art {I}.
\newblock {\em Journal of the IGPL}, 6:73--146, 1998.

\bibitem{gabbay2009quantification}
D.M. Gabbay, D.~Skvortsov, and V.~Shehtman.
\newblock {\em Quantification in nonclassical logic}.
\newblock Elsevier, 2009.

\bibitem{kracht1999tools}
M.~Kracht.
\newblock {\em Tools and techniques in modal logic}, volume 142.
\newblock Elsevier, 1999.

\bibitem{kremer2014quantified}
Ph. Kremer.
\newblock Quantified modal logic on the rational line.
\newblock {\em The Review of Symbolic Logic}, 7(3):439--454, 2014.

\bibitem{kudinov_aiml12}
A.~Kudinov.
\newblock Modal logic of some products of neighborhood frames.
\newblock 9:386--394, 2012.

\bibitem{kudinov2014neighbourhood}
A.~Kudinov.
\newblock Neighbourhood frame product {KxK}.
\newblock {\em Advances in Modal Logic}, 10:373--386, 2014.

\bibitem{Kudinov_KNproduct17}
A.~Kudinov.
\newblock Neighborhood-{K}ripke product of modal logics.
\newblock {\em Topology, Algebra, and Categories in Logic (TACL'17)}, 2017.

\bibitem{pacuit2017neighborhood}
E.~Pacuit.
\newblock {\em Neighborhood semantics for modal logic.}
\newblock Springer, 2017.

\bibitem{rasiowa1963metamathematics}
H.~Rasiowa and R.~Sikorski.
\newblock {\em Metamathematics of Mathematics}.
\newblock Polish Scientific Publishers, Warszawa, 1963.

\bibitem{Segerberg1971}
K.~Segerberg.
\newblock {\em An essay in classical modal logic}.
\newblock Filosofiska f\"oreningen och Filosofiska institutionen vid Uppsala
  universitet (Uppsala), 1971.

\end{thebibliography}

\vspace*{10pt}

%\address{Institute for Information Transmission Problems RAS\\
%	\hspace*{9pt}Bolshoy Karetny per. 19, build.1\\
%	\hspace*{18pt}Moscow 127051 Russia\\
%	{\it E-mail}: kudinov@iitp.ru\\[8pt]
%}

\clearpage
%\Appendix

\end{document}